\documentclass[reqno,centertags,12pt]{article} \usepackage{amsthm}
\usepackage{latexsym} \usepackage{eucal}
\usepackage{latexsym} \usepackage{eucal} \usepackage{srcltx}
\usepackage{amssymb,amsmath,amsfonts,amssymb}
\usepackage{amsmath,amsthm,amscd,amssymb}
\usepackage{latexsym}
\usepackage {verbatim}

\usepackage{graphics}
\usepackage{amsmath}
\usepackage{amsfonts,amssymb,amsthm,eucal}
\usepackage{pstricks,multido,graphicx,psfrag}
\usepackage{overpic}
\usepackage{amssymb,amsmath,amsfonts,amssymb}
-\marginparwidth \oddsidemargin -4mm \evensidemargin
\oddsidemargin
\advance\hoffset by 5mm

\def\@abssec#1{\vspace{.05in}\footnotesize \parindent .2in
{\bf #1. }\ignorespaces}

\def\1{{\bf 1}}

\chardef\set=35

\def\x{{\bf x}}

\def\u{{\bf u}}
\def\p{{\bf p}}

\def\n{{\bf n}}

\def\k{{\bf k}}

\def\s{{\bf s}}
\def\p{{\bf p}}

\def\m{{\bf m}}

\def\R{{\mathbb R}}
\def\Z{{\mathbb Z}}

\def\E{{\cal E}}

\def\embb{\supset \kern-11.6pt\lower.335ex\hbox{$\scriptscriptstyle <$} \;}

\title{Extended States for Polyharmonic Operators with Quasi-periodic Potentials in Dimension Two}
\author{Yu. Karpeshina, R. Shterenberg}

\date{}
\begin{document}

\maketitle

\begin{abstract} We consider a polyharmonic operator $H=(-\Delta)^l+V(\x)$ in
dimension two with $l\geq 2$, $l$ being an integer, and a
quasi-periodic potential $V(\x)$. We prove that the spectrum of $H$
contains a semiaxis and there is a family of generalized
eigenfunctions at every point of this semiaxis with the following
properties. First, the eigenfunctions are close to plane waves
$e^{i\langle \k,\x\rangle }$ at the high energy region. Second, the
isoenergetic curves in the space of momenta $\k$ corresponding to
these eigenfunctions have a form of slightly distorted circles with
holes (Cantor type structure). A new method of multiscale analysis
in the momentum space is developed to prove these
results.\end{abstract} We study  an operator
    \begin{equation}
    H=(-\Delta)^l+V(\x) \label{main0}
    \end{equation}
    in dimension two, where $l$ is an integer, $l\geq 2$, $V(\x)$ is a quasi-periodic
potential being a trigonometric polynomial:
\begin{equation}\label{V}
V=\sum\limits_{\s_1,\s_2\in\Z^2,\,0<|\s_1|+|\s_2|\leq
Q}V_{\s_1,\s_2}e^{2\pi i\langle \s_1+\alpha \s_2,\x\rangle},\ \
1\leq Q<\infty . \end{equation} We assume  that the irrationality
measure $\mu$ of $\alpha$ is finite: $\mu<\infty$, or in other
words, that $\alpha$ is not a Liouville number \footnote{ Note,  that $\mu\geq 2$ for any irrational number $\alpha $.}.


 The one-dimensional situation $d=1$, $l=1$ is thoroughly
   investigated in discrete and continuum settings,
   see e.g. \cite{DiSi}--\cite{FK} and references there.
   It is known that a one-dimensional  quasi-periodic Schr\"{o}dinger operator demonstrates spectral and transport
   properties which
   are \underline{not} close to those of a periodic operator.
    The spectrum of the quasi-periodic
   operator is, as a rule, a Cantor set, while in the periodic case, it has a band structure.
   In the periodic case the spectrum is absolutely continuous, while in the
   quasi-periodic case, it
   can have any nature: absolutely continuous, singular continuous and  pure point.
   The transition between different types  of spectrum can happen even with a small change of
   a coefficient in a quasi-periodic operator \cite{J1}.
    The mechanism of the difference in
   spectral behavior between periodic and quasi-periodic cases can be explained by a
   phenomenon
    which is known as resonance tunneling in quantum mechanics. It is associated with small
    denominators appearing in formal series of perturbation theory. Since the spectrum of
    the one-dimensional Laplacian is thin  (multiplicity 2),  resonance
    tunneling can produce an effect strong enough to destroy
    the spectrum. If a potential is periodic, then  resonance tunneling  produces
    gaps in the spectrum near the points $\lambda _n=(\pi n/a)^2$, $n\in \Z$, $a$ being the
    period of the potential. If the potential is quasi-periodic, then it
    can be thought as a sort of combination of  infinite number of periodic potentials,
     each
    of them producing gaps near its own $\lambda _n$-s. Since the set of all
    $\lambda _n$-s can be
    dense, the number of points surrounded by gaps can be dense too. Thus, the spectrum gets a
    Cantor like structure. The properties of the operator in the high energy region for the continuum case $d=1$ are studied in
   \cite{DiSi}-\cite{MP}, \cite{E}. The KAM method is used to prove absolute continuity of the
      spectrum and existence of quasiperiodic solutions at high
      energies.

    There are  important results on the density of states, spectrum,
     localization concerning the  quasi-periodic operators in $\Z^d$ and, partially, in
     $\R^d$, $d>1$, e.g. \cite{Sh1}--\cite{68a}.
        However, it is still much less known  about (\ref{main0}) then about its
    one-dimensional analog. The properties of the spectrum in the high energy region, existence of extended states and
    quantum transport are  still  wide open problems
    in the multidimensional case.

 Here we study properties of the spectrum and
eigenfunctions of (\ref{main0}) in the high energy region. We prove
the following results for the case $d=2$, $l\geq 2$.
    \begin{enumerate}
    \item The spectrum of the operator (\ref{main0})
    contains a semiaxis.

    This is a generalization of  a renown Bethe-Sommerfeld conjecture, which states that in the case of a periodic
    potential, $l=1$ and $d\geq 2$, the spectrum of \eqref{main0} contains a semiaxis.
    There is a variety of proofs for the periodic case, the earliest one is \cite{14r}. For a limit-periodic periodic potential, being periodic in one direction,  the conjecture is proved in \cite{SS}. For a general case
    of limit-periodic potential the conjecture is proven in \cite{KL1}-\cite{KL3}. Here we present the first proof of (a generalized) Bethe-Sommerfeld conjecture for a quasi-periodic potential.

    \item There are generalized eigenfunctions $\Psi_{\infty }(\k, \x )$,
    corresponding to the semi-axis, which are close to plane waves:
    for every $\k $ in an extensive subset $\cal{G} _{\infty }$ of
$\R^2$, there is
    a solution $\Psi_{\infty }(\k, \x)$ of the  equation
    $H\Psi _{\infty }=\lambda _{\infty }\Psi _{\infty }$ which can be
described by
    the formula:
    \begin{equation}
    \Psi_{\infty }(\k, \x)
    =e^{i\langle \k, \x \rangle}\left(1+u_{\infty}(\k,
    \x)\right), \label{qplane}
    \end{equation}
    \begin{equation}
    \|u_{\infty}\|_{L_{\infty }(\R^2)}\underset{|\k| \rightarrow
     \infty}{=}O\left(|\k|^{-\gamma _1}\right),\ \ \ \gamma _1>0,
    \label{qplane1}
    \end{equation}
    where $u_{\infty}(\k, \x)$ is a quasi-periodic
    function, namely a point-wise convergent series of exponentials $e^{i\langle\n+\alpha \m,\x\rangle}$, $\n ,\m \in \Z^2$.
    The  eigenvalue $\lambda _{\infty }(\k)$, corresponding to
    $\Psi_{\infty }(\k, \x)$, is close to $|\k|^{2l}$:
    \begin{equation}
    \lambda _{\infty }(\k)\underset{|\k| \rightarrow
     \infty}{=}|\k|^{2l}+
    O\left(|\k|^{-\gamma _2}\right),\ \ \ \gamma _2>0. \label{16a}
    \end{equation}
     The ``non-resonant" set $\cal{G} _{\infty }$ of
       vectors $\k$, for which (\ref{qplane}) -- (\ref{16a}) hold, is
       an extensive
       Cantor type set: ${\cal G} _{\infty }=\cap _{n=1}^{\infty }{\cal G}
_n$,
       where $\{{\cal G} _n\}_{n=1}^{\infty}$ is a decreasing sequence of
sets in $\R^2$. Each ${\cal G} _n$ has a finite number of holes  in each
bounded
       region. More and more holes appear when $n$ increases,
       however
       holes added at each step are of smaller and smaller size.
       The set $\cal{G} _{\infty }$ satisfies the estimate:
       \begin{equation}\left|\cal{G} _{\infty }\cap
        \bf B_R\right|\underset{R \rightarrow
     \infty}{=}|{\bf B_R}| \bigl(1+O(R^{-\gamma _3})\bigr),\ \ \ \gamma
_3>0,\label{full}
       \end{equation}
       where $\bf B_R$ is the disk of radius $R$ centered at the
       origin, $|\cdot |$ is the Lebesgue measure in $\R^2$.

       \item The set $\cal{D}_{\infty}(\lambda)$,
defined as a level (isoenergetic) set for $\lambda _{\infty }(
\k)$,
$$ {\cal D} _{\infty}(\lambda)=\left\{ \k \in \cal{G} _{\infty }
:\lambda _{\infty }(\k)=\lambda \right\},$$ is proven
to be a slightly distorted circle with infinite number of holes. It
can be described by  the formula:
 \begin{equation} {\cal
D}_{\infty}(\lambda)=\left\{\k:\k=\varkappa
_{\infty}(\lambda, \vec \nu){\vec \nu},
    \ {\vec \nu} \in {\cal B}_{\infty}(\lambda)\right
    \}, \label{Dinfty}
    \end{equation}
where ${\cal B}_{\infty }(\lambda )$ is a subset of the unit circle
$S_1$. The set ${\cal B}_{\infty }(\lambda )$ can be interpreted as
the set of possible  directions of propagation for  almost plane
waves (\ref{qplane}). The set ${\cal B}_{\infty }(\lambda )$ has a
Cantor type structure and an asymptotically full measure on $S_1$ as
$\lambda \to \infty $:
\begin{equation}
L\bigl({\cal B}_{\infty }(\lambda )\bigr)\underset{\lambda
\rightarrow
     \infty}{=}2\pi +O\left(\lambda^{-\gamma _4/2l}\right),\ \ \
     \gamma_4>0,
\label{B}
\end{equation}
here and below $L(\cdot)$ is a length of a curve. The value
$\varkappa _{\infty }(\lambda ,{\vec \nu} )$ in (\ref{Dinfty}) is the
``radius" of ${\cal D}_{\infty}(\lambda)$ in a direction ${\vec \nu} $.
The function $\varkappa _{\infty }(\lambda ,{\vec \nu} )-\lambda^{1/2l}$
describes the deviation of ${\cal D}_{\infty}(\lambda)$ from the
perfect circle of the radius $\lambda^{1/2l}$. It is proven that the
deviation is asymptotically small:
\begin{equation}
\varkappa _{\infty }(\lambda ,{\vec \nu} )\underset{\lambda
\rightarrow
     \infty}{=}\lambda^{1/2l}+O\left(\lambda^{-\gamma _5 }\right),
\ \ \ \gamma _5>0. \label{h}
\end{equation}
\item The branch of the
    spectrum  corresponding to $\Psi_{\infty }(\k,
\x)$ (the semiaxis) is absolutely continuous.

\end{enumerate}

To prove the results listed above we  suggest a method which can be
described as {\em multiscale analysis in the space of momenta}. This
is a development of the method, which is used in \cite{KL1}--\cite{KL3} for the
case of limit-periodic potentials. The essential difference is that
in \cite{KL1}--\cite{KL3} we constructed a modification of KAM method, where the
space variable  $\x$ still plays some role (e.g. in the uniform  in
$\x$ approximation of a limit-periodic potential by periodic ones),
while in the present situation all considerations are happening in
the space of the dual variable $\k$. The KAM method in \cite{KL1}--\cite{KL3}
is  motivated by \cite{2}--\cite{3}, where the method is used for
periodic problems. Multiscale analisys which we apply here is deeply
analogous to the original multiscale method developed in \cite{FrSp}
(see also \cite{BG}, \cite{74}) for the proof of localization. The
essential difference is that
 in \cite{FrSp}, \cite{BG}, \cite{74} the multiscale procedure is constructed with respect to space variable $\x$ to prove localization, while we construct a multiscale procedure in the space of momenta $\k$ to prove delocalization.

        Here is a brief description of the iteration procedure which leads to the results described above.
        Indeed, let $\k \in \R^2$. We consider a set of finite linear combinations of plane waves $e^{i\langle\k+\p+\alpha \m,\x\rangle}$,
        $ \p,\m \in \Z^2$.
         The set is invariant under action of the differential expression (\ref{main0}). Let $H(\k)$ be a
         matrix describing action of  (\ref{main0}) in the linear set of the exponentials.
          Obviously,
          $$H(\k)=H_0(\k)+V,\ \  H_0(\k)_{(\p,\m), (\p',\m ')}=|\k+\p+\alpha \m|_{\R^2}^{2l}\delta _{(\p,\p')}
          \delta_{ (\m,\m')},$$ $$V_{(\p,\m), (\p',\m')}=V_{\p-\p', \m-\m'}.$$
          Next, we consider  an expanding sequence of finite sets $M_n$ in the space $\Z^2\times \Z^2$ of indices
          $(\p,\m)$: $M_n\subset M_{n+1}$, $\lim _{n\to \infty } M_n= \Z^2\times \Z^2$. Let $P_n$ be
          the characteristic projection of  set $M_n$ in the space $\ell^{2}(\Z^2\times \Z^2)$. We consider a sequence
          of finite matrices $H^{(n)}(\k)=P_nH(\k)P_n$. Each matrix corresponds to a finite dimensional operator in
          $\ell^{2}(\Z^2\times \Z^2)$, given that the operator acts as zero on $(I-P_n){\ell}^2$.
          For each $n$ we construct a ``non-resonant" set ${\cal G}_n$ in the space $\R^2$ of momenta $\k$, such that:
          if $\k \in {\cal G}_n$, then $H ^{(n)}(\k)=P_nH(\k)P_n$ has an eigenvalue $\lambda _n(\k)$ and
          its spectral projector $\E _n(\k)$ which can be described by perturbation formulas with respect to the previous
          operator $H ^{(n-1)}(\k)$. If $\k \in \cap _{n=1}^{\infty} {\cal G} _n$ then $\lambda _n(\k)$ and $\E _n(\k)$
          have  limits. The linear combinations of the exponentials, corresponding to the projectors $\E _n(\k)$, have a
          point-wise limit in $\x$, the limit being a generalized eigenfunction of \eqref{main0}.
         The generalized eigenfunction is close to the plane wave $e^{i<\k,\x>}$ in the high energy region.

          Each matrix $H ^{(n)}$ is considered as a perturbation of a matrix $\hat H ^{(n)}$, the latter has a block
          structure, i.e., consists of a variety of blocks $H^{(s)}(\k +\p +\alpha \m)$, $s=1,...,n-1$, and, naturally,
          some diagonal terms. Blocks with different indices $(s)$ have    sizes of different orders of magnitude  (the size increasing
          with $s$). Thus we have a multiscale structure  in the definition of $\hat H ^{(n)}$. We use $\hat H ^{(n)}(\k)$ as a starting operator to construct perturbation series for $H^{(n)}(\k)$.   At a step $n$ we apply our knowledge of spectral properties of $H^{(s)}(\k +\p'+\alpha \m')$, $s=1,...,n-1$, $\p', \m' \in  \Z^2$, obtained in the previous steps, to
        describe spectral properties of  $H^{(n)}(\k +\p +\alpha \m)$, $\p,\m\in \Z^2$ and to construct ${\cal G} _n$.

At step one we use a regular perturbation theory and elementary geometric considerations to prove the following results. There is a set ${\cal G}_1\subset \R^2$ such that: if $\k\in {\cal G}_1$, then the operator $H^{(1)}(\k)$ has a single eigenvalue close to the unperturbed one:
\begin{equation} \lambda ^{(1)}(\k)\underset{|\k| \rightarrow
     \infty}{=}|\k|^{2l}+
    O\left(|\k|^{-\gamma _2}\right),\ \ \ \gamma _2>0.\label{lambda-1} \end{equation}
     A normalized eigenvector $\u^{(1)}$ is also close to the unperturbed one: $\u^{(1)}=\u^{(0)}+\tilde \u^{(1)}$, where
     $(\u^{(0)})_{(\p,\m)}=\delta _{\p,{\bf 0}}\delta_{ \m,{\bf 0}}$ and the $l^{1}$-norm of $\tilde \u^{(1)}$ is small:
     $\|\tilde \u^{(1)}\|_{l^{1}}<|\k|^{-\gamma _1}$, $\gamma _1>0$. It follows that:
 \begin{equation}
    \Psi_1 (\k, \x)
    =e^{i\langle \k, \x \rangle}+\tilde u_{1}(\k,
    \x),\ \ \
\|\tilde u_{1}\|_{L_{\infty }(\R^2)}\underset{|\k| \rightarrow
     \infty}{=}O(|\k|^{-\gamma _1}),\ \
\ \gamma _1>0, \label{na}
\end{equation}
    where  $\Psi_1 (\k, \x)$, $\tilde u_{1}(\k, \x)$ are the linear combinations of the exponentials corresponding to
    vectors $\u^{(1)}$ and $\tilde \u^{(1)}$, respectively. It is shown that function $\Psi_1 (\k, \x)$
    satisfies the equation for eigenfunctions with a good accuracy:
    \begin{equation} \label{eqforeigenfunctions-1}-\Delta \Psi_1+V\Psi_1=|\k|^{2l}\Psi_1+f_1,\ \ \|f_1\|_{L_{\infty }(\R^2)}\underset{|\k| \rightarrow
     \infty}{=}O(|\k|^{-\gamma _6})\ \
\ \gamma _6>0.
\end{equation}
Relation \eqref{lambda-1} is differentiable:
\begin{equation} \nabla \lambda ^{(1)}(\k)\underset{|\k| \rightarrow
     \infty}{=}2l|\k|^{2l-2}\k+
    O\left( |\k|^{-\gamma _7}\right),\ \ \ \gamma _7>0.\label{dlambda-1} \end{equation}
     Next, we construct a sequence ${\cal G}_n$, $n\geq 2$, such  for any $\k \in {\cal G}_n$ the operator $H^{(n)}(\k)$ has a single eigenvalue $\lambda ^{(n)}(\k)$ in a super exponentially small neighborhood of $\lambda ^{(n-1)}(\k)$:
\begin{equation} \lambda ^{(n)}(\k)\underset{|\k| \rightarrow
     \infty}{=}\lambda ^{(n-1)}(\k)
   + O\left(|\k|^{-|\k|^{\gamma _8n}}\right),\ \ \ \gamma _8>0.\label{lambda-n} \end{equation}
   Similar estimates hold for the eigenvectors and the  corresponding functions $\Psi_n(\k,\x)$:
   \begin{equation}
   \Psi _n(\k,\x)=\Psi _{n-1}(\k,\x)+\tilde u _n(\k,\x),\ \ \ \|\tilde u_{n}\|_{L_{\infty }(\R^2)}\underset{|\k| \rightarrow
     \infty}{=} O\left(|\k|^{-|\k|^{\gamma _9n}}\right),\ \ \ \gamma _9>0.\ \ \label{na-n}
\end{equation}
\begin{equation} \label{eqforeigenfunctions-1}-\Delta \Psi_n+V\Psi_n=\lambda ^{(n)}(\k)\Psi_n+f_n,\ \ \|f_n\|_{L_{\infty }(\R^2)}\underset{|\k| \rightarrow
     \infty}{=}O\left(|\k|^{-|\k|^{\gamma _{10}n}}\right),\ \ \
\ \gamma _{10}>0.
\end{equation}
Formula \eqref{lambda-n} is differentiable with respect to $\k$:
\begin{equation} \nabla \lambda ^{(n)}(\k)\underset{|\k| \rightarrow
     \infty}{=}\nabla \lambda ^{(n-1)}(\k)+O\left(|\k|^{-|\k|^{\gamma _8n}}\right),\ \ \
\ \gamma _8>0.\label{dlambda-n} \end{equation}
In fact, for large $n$ estimates \eqref{lambda-n} -- \eqref{dlambda-n} are even stronger.
\begin{figure}
\begin{minipage}[t]{8cm}
\centering
\includegraphics
[totalheight=.25\textheight]{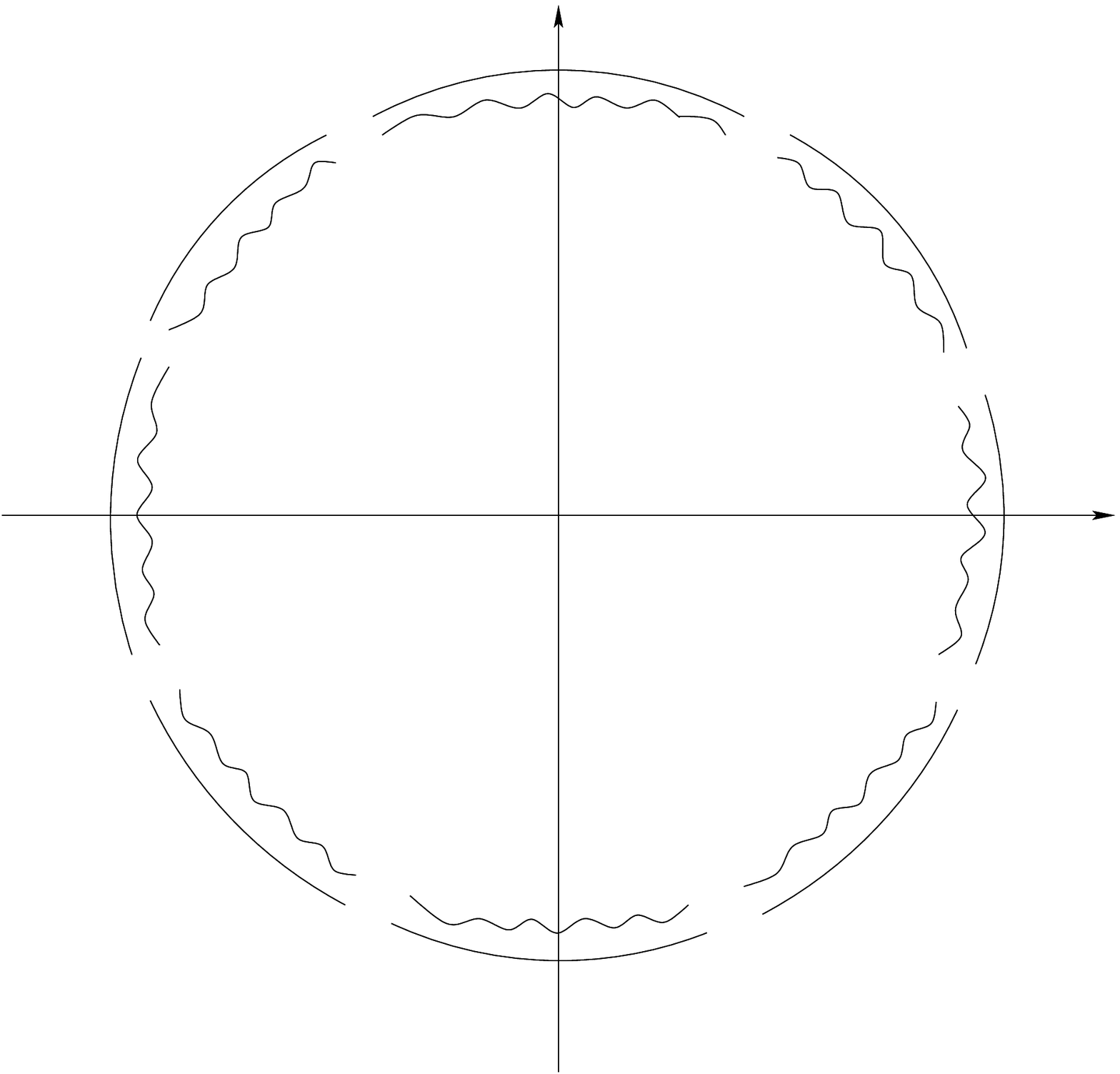} \caption{Isoenergetic curve
${\cal D}_1(\lambda)$}\label{F:1}
\end{minipage}
\hfill
\begin{minipage}[t]{8cm}
\centering
\includegraphics[totalheight=.25\textheight]{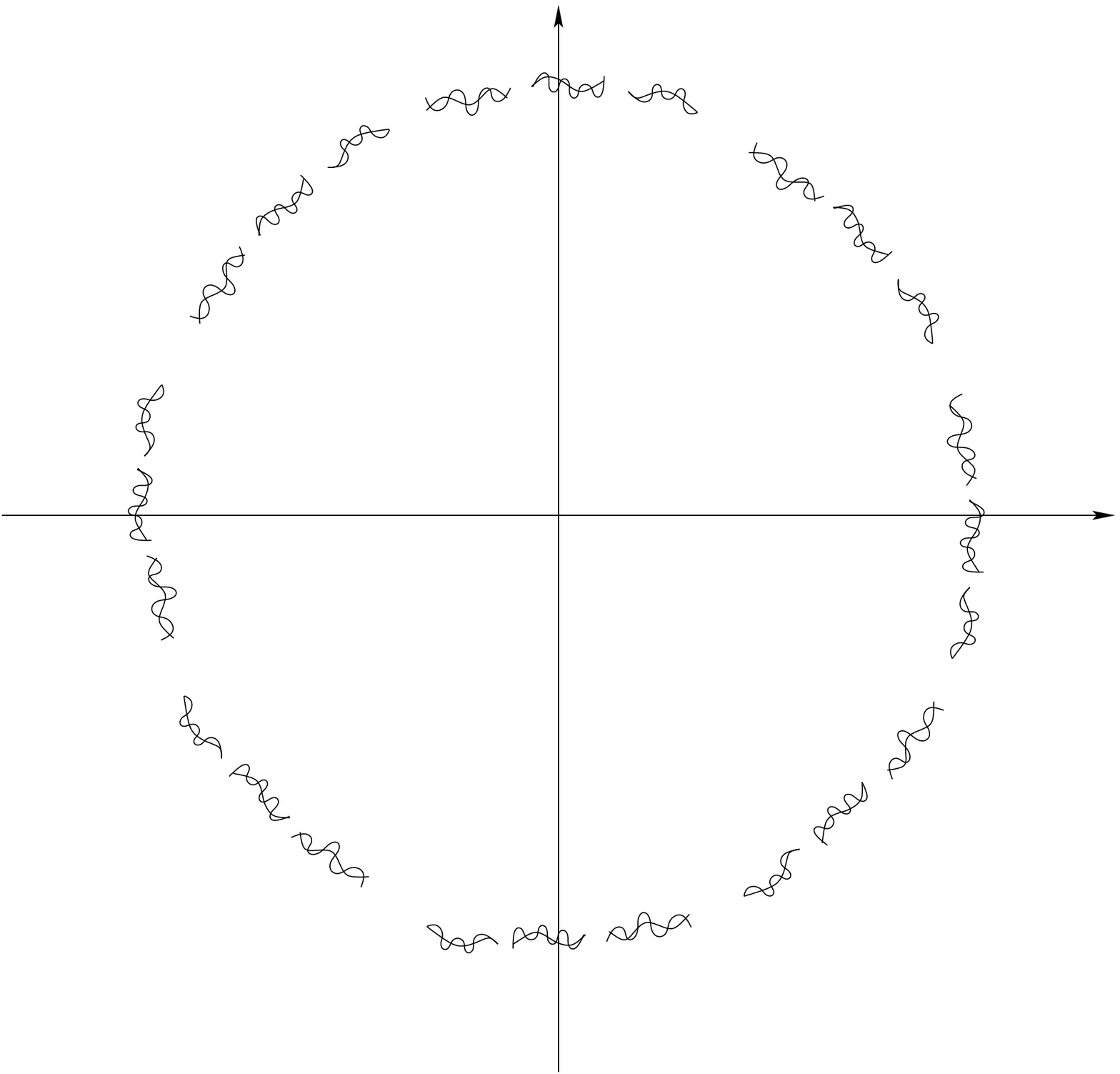}
\caption{Isoenergetic curve ${\cal D}_2(\lambda)$}\label{F:2}
\end{minipage}\hfill
\end{figure}
The non-resonant set ${\cal G} _{n}$
        is proven to be
       extensive in $\R^2$:
       \begin{equation}
       \left|{\cal G} _{n}\cap
       \bf B_R\right|\underset{R \rightarrow
     \infty}{=}|{\bf B_R}|\bigl(1+O(R^{-\gamma _3})\bigr). \label{16b}
       \end{equation}
       Estimates (\ref{lambda-n}) -- (\ref{16b}) are uniform in $n$.

The set ${\cal D}_{n}(\lambda)$ is defined as the level
(isoenergetic) set for the non-resonant eigenvalue $\lambda ^{(n)}(\k)$:
$$ {\cal D} _{n}(\lambda)=\left\{ \k \in {\cal G} _n:\lambda ^{(n)}(\k)=\lambda \right\}.$$
This set is proven to be a slightly distorted circle with a
finite number of holes (see Fig. \ref{F:1}, \ref{F:2}). The set
${\cal D} _{n}(\lambda)$ can be described by the formula:
\begin{equation}
{\cal D}_{n}(\lambda)=\left\{\k:\k=
    \varkappa_{n}(\lambda, {\vec \nu})\vec \nu ,
    \ \vec \nu  \in {\cal B}_{n}(\lambda)\right\}, \label{Dn}
    \end{equation}
where ${\cal B}_{n}(\lambda )$ is a subset  of the unit circle
$S_1$. The set ${\cal B}_{n}(\lambda )$ can be interpreted as the
set of possible directions of propagation for  almost plane waves $\Psi _n(\k,\x)$, see (\ref{na}),
(\ref{na-n}). It has an asymptotically full measure on $S_1$ as
$\lambda \to \infty $:
\begin{equation}
L\bigl({\cal B}_{n}(\lambda )\bigr)\underset{\lambda \to \infty
}{=}2\pi +O\left(\lambda^{-\gamma _4/2l}\right). \label{Bn}
\end{equation}
Each set ${\cal
B}_{n}(\lambda)$ has only a finite number of holes, however their
number is growing with $n$. More and more holes of a smaller and
smaller size are added at each step. The value
$\varkappa_{n}(\lambda ,{\vec \nu} )-\lambda^{1/2l}$ gives the
deviation of ${\cal D}_{n}(\lambda)$ from the perfect circle of
the radius $\lambda^{1/2l}$  in the direction ${\vec \nu} $. It is
proven that the deviation is asymptotically small:
\begin{equation}
\varkappa_{n}(\lambda ,{\vec \nu}) =\lambda^{1/2l}+O\left(\lambda^{-
\gamma _5}\right),\ \ \ \ \frac{\partial \varkappa_{n}(\lambda
,{\vec \nu})}{\partial \varphi }=O\left(\lambda^{- \gamma _{11} }\right),\
\ \gamma_5,\gamma _{11}>0,\label{hn}
\end{equation}
$\varphi $ being an angle variable, ${\vec \nu} =(\cos \varphi ,\sin
\varphi )$.  Estimates (\ref{Bn}), (\ref{hn}) are uniform in $n$.

On each step more and more points are excluded from the
non-resonant sets ${\cal G} _n$, thus $\{ {\cal G} _n \}_{n=1}^{\infty }$ is a
decreasing sequence of sets. The set ${\cal G} _\infty $ is defined as the
limit set: ${\cal G} _\infty=\cap _{n=1}^{\infty }{\cal G} _n $. It
has an infinite number of holes, but nevertheless satisfies the
relation (\ref{full}). For every $
\k \in {\cal G} _\infty $ and every $n$, there is a generalized
eigenfunction of $H^{(n)}$ of the type  (\ref{na}), (\ref{na-n}). It is
proven that  the sequence of
$\Psi _n(\k, \x)$ has a limit in $L_{\infty }(\R^2)$ when $
\k \in {\cal G} _\infty $.
The function $\Psi _{\infty }(\k, \x)
=\lim _{n\to \infty }\Psi _n(\k, \x)$ is a generalized
eigenfunction of $H$. It can be written in the form
(\ref{qplane}) -- (\ref{qplane1}).
Naturally, the corresponding eigenvalue $\lambda _{\infty }(\k) $ is
the limit of $\lambda ^{(n)}(\k )$ as $n \to \infty $.

It is shown that $\{{\cal B}_n(\lambda)\}_{n=1}^{\infty }$  is a
decreasing sequence of sets,  on each step more and more directions
being excluded. We consider the limit ${\cal B}_{\infty}(\lambda)$
of ${\cal B}_n(\lambda)$:
    \begin{equation}{\cal B}_{\infty}(\lambda)=\bigcap_{n=1}^{\infty} {\cal
    B}_n(\lambda).\label{Dec8a}
    \end{equation}
    This set has a Cantor type structure on the unit circle.
    It is proven that ${\cal B}_{\infty}(\lambda)$ has an asymptotically
    full measure on the unit circle (see (\ref{B})).
    We prove
    that the sequence $\varkappa _n(\lambda ,{\vec \nu} )$, $n=1,2,... $,
describing the
     isoenergetic curves ${\cal D}_n(\lambda)$, quickly converges as $n\to
\infty$. We show that ${\cal D}_{\infty}(\lambda)$ can be described as the
limit of  ${\cal D}_n(\lambda)$ in the sense (\ref{Dinfty}), where
$\varkappa _{\infty}(\lambda, \vec \nu )=\lim _{n \to \infty} \varkappa
_n(\lambda, \vec \nu )$ for every $\vec \nu  \in {\cal B}_{\infty}(\lambda)$.
It is shown that the derivatives of the functions $\varkappa
_n(\lambda, \vec \nu )$ (with respect to the angle variable on the unit
circle) have a limit as $n\to \infty $ for every $\vec \nu  \in {\cal
B}_{\infty}(\lambda)$. We denote this limit by $\frac{\partial
\varkappa_{\infty}(\lambda ,\vec \nu)}{\partial \varphi }$. Using
(\ref{hn}), we  prove that
    \begin{equation}\frac{\partial \varkappa_{\infty}(\lambda ,\vec \nu)}{\partial
\varphi }=O\left(\lambda^{- \gamma _{11} }\right).\label{Dec9a}
\end{equation} Thus, the limit curve ${\cal D}_{\infty}(\lambda)$ has a
tangent vector in spite of its Cantor type structure, the tangent
vector being the limit of corresponding tangent vectors for ${\cal
D}_n(\lambda)$ as $n\to \infty $.  The curve  ${\cal
D}_{\infty}(\lambda)$ looks as
  a slightly distorted circle with
infinite number of holes for every sufficiently large $\lambda $, $\lambda >\lambda _*(V)$. It immediately follows that $[\lambda _*, \infty )$ is in the spectrum of $H$ (Bethe-Sommerfeld conjecture).

The main technical difficulty to overcome is the construction of
   non-resonant sets
${\cal B} _n(\lambda)$ for every fixed sufficiently large $\lambda $,
$\lambda >
\lambda_0 (V)$,
where $\lambda_0(V)$ is the same for all $n$. The set
${\cal B} _n(\lambda)$ is obtained  by deleting a ``resonant" part
from
${\cal B}_{n-1}(\lambda)$.
Definition of  ${\cal B} _{n-1}(\lambda)\setminus {\cal B}_{n}(\lambda)$
includes
eigenvalues of $H^{(n-1)}(\k)$. To describe  ${\cal B} _{n-1}(\lambda)\setminus
{\cal B}_{n}(\lambda)$ one
has to consider
  not only non-resonant
eigenvalues of the type (\ref{lambda-1}),  (\ref{lambda-n}), but also
resonant
eigenvalues, for which no suitable  formulas are known. Absence of
formulas causes difficulties in estimating the size of ${\cal B}
_{n-1}(\lambda)\setminus {\cal B}_n(\lambda)$.
       To treat this problem we start with introducing
        an angle variable $\varphi \in [0,2\pi
       )$,  $\vec \nu  = (\cos \varphi ,
       \sin \varphi )\in S_1$ and consider sets ${\cal B}_n(\lambda)$ in
terms of this
       variable.
        Next, we show that
the resonant set  ${\cal B} _{n-1}(\lambda)\setminus
{\cal B}_{n}(\lambda)$ can be described as the set of zeros of  functions of the type
        $$ \det \Bigl(H^{(s)}\bigl(\vec \varkappa _{n-1}(\varphi )+\p+\alpha \m
       \bigr)-\lambda-\varepsilon \Bigr), \ \ \ s=1,...,n-1,\ \ \ (\p,\m)\in M_n\setminus{(\bf 0,\bf 0)},$$
       where $\vec \varkappa _{n-1}(\varphi )$ is a vector-function
describing ${\cal D}_{n-1}
       (\lambda)$: $\vec \varkappa _{n-1}(\varphi )=\varkappa _{n-1}(\lambda
,{\vec \nu} ){\vec \nu} $. To obtain ${\cal B} _{n-1}(\lambda)\setminus
       {\cal B}_{n}(\lambda)$ we take all values of $\varepsilon $ in a small
interval
       and $(\p,\m)$ in some subset of $M_n$.
        Further, we extend our
       considerations to
        a complex neighborhood $\varPhi _0$ of $[0,2\pi
       )$. We show that the    determinants are analytic functions of
        $\varphi $ and, by this,
         reduce the
        problem of estimating the size of the resonant set
         to
        a problem in complex analysis. We use theorems for analytic
functions to count
          zeros of the determinants and to investigate how far  the zeros move
when
         $\varepsilon $ changes. It enables us to estimate the size
         of the zero set of the determinants, and, hence, the size of
         the non-resonant set $\varPhi _n\subset \varPhi _0$,  which is
defined as a
         non-zero set for the determinants.
          Proving that the non-resonant set $\varPhi _n$
         is sufficiently large, we
          obtain estimates (\ref{16b}) for ${\cal G} _n$ and (\ref{Bn}) for
         ${\cal B}_n$, the set  ${\cal B}_n$ corresponding to the real part of
         $\varPhi _n$.
 \begin{figure}
\begin{minipage}[t]{8cm}
\centering \psfrag{Phi_2}{$\Phi_2$}
\includegraphics
[totalheight=.25\textheight]{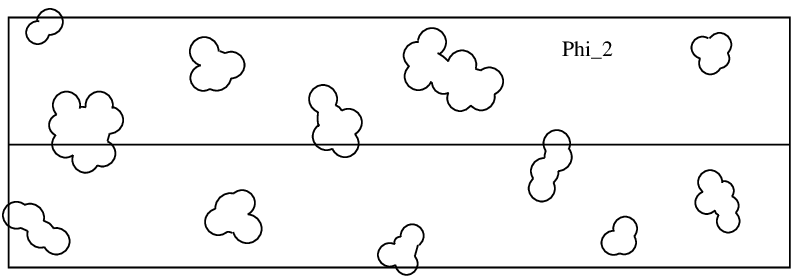} \caption{Set $\varPhi
_2$}\label{F:3}
\end{minipage}
\hfill
\end{figure}
\begin{figure}
\begin{minipage}[t]{8cm}
\centering \psfrag{Phi_3}{$\Phi_3$}
\includegraphics[totalheight=.25\textheight]{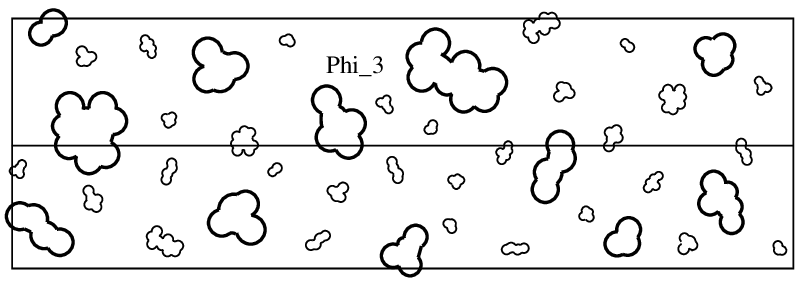}
\caption{Set $\varPhi _3$}\label{F:4}
\end{minipage}\hfill
\end{figure}
To obtain $\varPhi _n$ we delete  from $\varPhi _0$ more and more discs
(holes) of smaller and
         smaller radii at each step. Thus, the non-resonant set $\varPhi
         _n\subset \varPhi _0$ has a structure of Swiss
         Cheese (Fig. \ref{F:3}, \ref{F:4}). Deleting  a resonance set from $\varPhi _0$ at each
         step of the recurrent procedure  we call a ``Swiss
         Cheese Method".  The essential
         difference of our method from constructions of non-resonant sets  in similar
         situations before (see e.g. \cite{2}--\cite{3}, \cite{B3}) is that
we construct a
         non-resonant set not only in the whole space of a parameter
         ($\k\in \R^2$ here), but also on   isoenergetic curves
         ${\cal D}_n(\lambda )$ in
         the space of the parameter, when $\lambda $ is sufficiently large.
Estimates for the
         size of non-resonant sets on a curve require more subtle
         technical considerations than those sufficient for
         description of a non-resonant set in the whole space of
         the parameter. But as a reward, such estimates enable us to
         show that every isoenergetic set for $\lambda>\lambda_0$ is not empty and
         thus, to prove Bethe-Sommerfeld conjecture.

  Note that generalization of  the results from the case $l>1$, $l$ being an integer, to the case of rational $l$ satisfying the
         same inequality is relatively simple; it requires just slightly more careful technical considerations. The restriction $l>1
         $  is also technical, though it is more difficult to lift. The condition $l>1$ is needed only for
         the second step of the recurrent procedure. The authors plan to consider
         the case $l=1$ in a forthcoming paper.  The
         requirement $\mu < \infty $ is  essential, since we use it to estimate the minimal values of $|\p+\alpha \m|$ when $(\p, \m)\in M_n\setminus{(\bf 0,\bf 0)}$. Such estimates are necessary for controlling small denominators in the perturbation series at each step.

\vspace{5mm} \noindent {\bf Acknowledgement} The authors are very
grateful to Prof. Leonid Parnovski for useful discussions and to
Prof. Young-Ran Lee for allowing us to use figures 1-4 from
\cite{KL1}.

\vskip 0.3in

Department of Mathematics, University of Alabama at Birmingham,

1300 University Boulevard, Birmingham, AL 35294.

e-mail: karpeshi@math.uab.edu; shterenb@math.uab.edu

\end{document}